\documentclass[a4paper]{article}

\author{Joaquim Ro\'e
   \thanks {Partially supported by CIRIT 1997FI-00141,
         CAICYT PB95-0274,
         and ``AGE-Algebraic Geometry in Europe" contract no. ERB940557.
         } \\
  \small{Departament d'\`Algebra i Geometria, 
           Universitat de Barcelona,} \\
  \small{Gran Via, 585, E-08007, Barcelona. }\\
  \small{e-mail: jroevell@cerber.mat.ub.es} \\
           }

\usepackage{amstex}
\usepackage{amsthm}
\usepackage{amssymb}
\usepackage{xypic}
\usepackage{graphics}

\newtheorem{Def}{Definition}[section]
\newtheorem{Lem}[Def]{Lemma}
\newtheorem{Cor}[Def]{Corollary}
\newtheorem{Pro}[Def]{Proposition}
\newtheorem{Teo}[Def]{Theorem}

\newcommand{\BP}{\operatorname{BP}}

\newcommand{\Hilb}{\operatorname{Hilb}}
\newcommand{\Spec}{\operatorname{Spec}}
\newcommand{\MaxSpec}{\operatorname{MaxSpec}}
\newcommand{\len}{\operatorname{length}}

\renewcommand{\P}{\mathbb{P}}
\renewcommand{\H}{{\cal H}}
\renewcommand{\O}{{\cal O}}
\renewcommand{\L}{{\cal L}}
\newcommand{\V}{{\cal V}}
\newcommand{\I}{{\cal I}}
\newcommand{\M}{\mathfrak m}
\newcommand{\C}{\mathbb{C}}
\newcommand{\Z}{\mathbb{Z}}
\newcommand{\G}{\mathbb{G}}
\newcommand{\m}{{\bf m}}
\renewcommand{\r}{{\bf r}}

\begin{document}         

\title{Tacnodes and cusps}
\maketitle

%
%
\section {Introduction}

The study of linear systems of algebraic plane curves with fixed
imposed singularities is a classical subject which has
recently experienced important progress. The Horace
method introduced by A. Hirschowitz in \cite{hir} has
been successfully exploited to prove many $\smash{H^1}$--vanishing
theorems, even in higher dimension. Other specialization
techniques, which include degenerations of the plane,
are due to Z. Ran \cite{ran1}, \cite{ran2} and C. Ciliberto and R. Miranda
\cite{cm2}. In \cite{gls}, G. M. Greuel, C. Lossen and
E. Shustin use a local specialization procedure together
with the Horace method to give the first asymptotically
proper general existence criterion for singular curves
of low degree. In this paper we develop a specialization
method which allows us to compute the dimension of several
linear systems as well as to substantially improve the
bounds of \cite{gls} for curves with tacnodes and cusps.

A closed subscheme $Z \subset \P^2$ is said to have
\emph{maximal rank in degree $d$} if the canonical map
$H^0(\P^2,\O_{\P^2}(d)) \rightarrow H^0(Z,\O_Z(d))$ has maximal
rank (cf. \cite{hir}). If $Z$ is a zero--dimensional scheme of length
$N$, this means that either there are no curves of
degree $d$ containing $Z$ or $Z$ imposes $N$ independent linear
conditions to curves of degree $d$. If $Z$ has maximal rank in
all degrees then we say simply that it has maximal rank.

For some classes of zero-dimensional schemes, it is known
that general members have maximal rank. For example, a general
union of double points has maximal rank. In many other cases,
however, a maximal rank statement has been conjectured only
(cf. \cite{bh2}, \cite{hir2}, \cite{mir}). We will
consider schemes
$Z=Z_1 \cup Z_2 \cup \ldots \cup Z_\rho$ where $Z_i$
are unibranched cluster schemes all whose points but one, whether
proper or infinitely near, are taken with multiplicity $\leq 2$.
In theorem \ref{rang} we prove that, under some mild numerical
conditions, a scheme $Z$ as above whose points are in general
position has maximal rank. This generalizes and unifies a range
of previously known results:
\begin{itemize}
\item In \cite{hir}, A. Hirschowitz proves that a
general union of distinct points with multiplicity 2 or 3 has maximal
rank. Our result generalizes the multiplicity 2 case by allowing
infinitely near points, which are not easily dealt with by the
Horace method.
\item M. V. Catalisano and A. Gimigliano in \cite{cat} and
C. Ciliberto and R. Miranda in \cite{cm1} prove that a union
of general unibranched cluster schemes all whose points have multiplicity
one (curvilinear schemes) has maximal rank. Our result generalizes this
by allowing multiplicity 2 points. 
\item In \cite{cm2}, C. Ciliberto and R. Miranda consider
quasi-homogeneous schemes, which consist of distinct multiple points,
all whose multiplicities but one are taken to be equal. They prove
a quite general maximal rank theorem which includes the computation
of the superabundant systems (corresponding to multiplicities for
which the general scheme does not have maximal rank). In the
particular case in which one has a point of multiplicity
not bigger than 5 and all other points have multiplicity 2,
we  extend their result by allowing infinitely near points.
\end{itemize}

Other results concerning linear systems with infinitely near base
points can be found in the literature, mainly in \cite{bh1},
\cite{bh3}, \cite{gls} and \cite{los}. Harbourne's results deal
with clusters of points lying on conics and cubics, which is not
necessary to us, and the $\smash{H^1}$--vanishing
of  \cite{gls} and \cite{los} is weaker than proving maximal rank.

The maximal rank theorem will allow us to prove the
existence of irreducible curves of low degree with
tacnodes and higher order cusps. The
reasoning is similar to that of \cite{bar}, \cite{gls}
or \cite{los}. Theorem \ref{existencia} is much sharper than
the result by Greuel, Lossen and Shustin (which nevertheless
applies to any kind of plane singularity) and also
sharper than the one Lossen obtains for tacnodes and
cusps (A-singularities), in part because our specialization
avoids the use of the Viro method to ``glue'' the singularities.
On the other hand, Barkats' result is a little bit sharper
than ours, but it is restricted to nodes and ordinary cusps
only.

\section {Preliminaries}
\label{preliminars}

Let $k$ be an algebraically closed field of characteristic zero,
$p$ a smooth point of a surface $S$ defined over $k$, $\O=\O_{S,p}$
the local ring of $p$ on $S$, $\M=\M_p$ its maximal ideal.

Consider a sequence of blowing-ups
$$
S_r \overset{\sigma_r}\longrightarrow S_{r-1} \longrightarrow \cdots
\overset{\sigma_2}\longrightarrow S_1 \overset{\sigma_1}\longrightarrow S_0=S
$$
where $\sigma_1$ is the blowing-up of $p$ and for $i>1$ the center
of $\sigma_i$ is a point $p_i$ which lies on the
exceptional divisor of $\sigma_{i-1}$. The sequence
$K=(p_1=p,p_2,\ldots,p_r)$ is a cluster with origin at $p$
for which every point is infinitely near the preceding one;
we will call these clusters \emph{unibranched}.
We write $S_K=S_r$ and $\pi_K:S_K \rightarrow S$
the composition of the blowing-ups. Usual facts
known for clusters hold in particular for unibranched
clusters; we now review some of them, referring the
reader to \cite{cma}, \cite{cll} for the proofs.

A point $p_j$ is said to be \emph{proximate} to $p_i$, $j>i$
if and only if it lies on the exceptional divisor
of blowing up $p_i$ (that is, $j=i+1$)
or on its \emph{strict} transform when $j>i+1$.
Every point in a unibranched cluster is proximate either
to one or to two points, except for $p_1$, which is
proximate to no one; if $p_j$ is proximate to $p_i$
then all the points between them are also proximate to $p_i$.
If $p_j$ is proximate to $p_i$ then there is a unique
point in $E_j$ proximate to both $p_i$ and $p_j$; otherwise
there is none.
If a point is proximate to two points, it is called
a \emph{satellite}, otherwise it is \emph{free}.
When considering more than one cluster at a time, we
will write $p_i(K)$ and $S_i(K)$ for the $i$--th point
of the cluster $K$ and the surface obtained by blowing
up the first $i$ points of $K$.

A \emph{system of multiplicities} for a cluster
$K=(p_1, p_2, \ldots, p_r)$ is a sequence
of integers $\m=(m_1, m_2, \ldots, m_r)$, and a pair
$(K,\m)$ where $K$ is a cluster and $\m$ a system
of multiplicities is called a \emph{weighted cluster}.
A system of multiplicities like
$$(m, 2, 2, \overset{\underset{\smile} i} \dots, 2,
     1, 1, \overset{\underset{\smile} j} \dots, 1,
     0, 0, \dots, 0)$$
will be designed as $(m,2^i,1^j)$.
Given a weighted cluster, there are an ideal and
a zero--dimensional subscheme of $S$ associated to it.
Let $E_i$ be the pullback (total transform) in $S_K$
of the exceptional divisor of blowing up $p_i$.
Then the ideal sheaf
$$ \H_{K,\m}=\left(\pi_K\right)_*
   \O_{S_K}(-m_1 E_1 - m_2 E_2 - \cdots - m_r E_r)$$
is supported at $p$, its stalk at $p$
is a complete $\M$--primary ideal $H_{K,\m} \subset \O$ 
and defines a zero--dimensional subscheme $Z_{K,\m}$ of $S$.
For $K$ unibranched we call $Z_{K,\m}$ a \emph{unibranched
cluster scheme}.
As an aside, note that if $I \subset \O$ is a complete
$\M$--primary ideal then there is
a weighted cluster $(K,\m)$ such that $I=H_{K,\m}$,
but this cluster does not need to be unibranched.

The same (unibranched) cluster scheme is sometimes defined
by different (unibranched) weighted clusters.
In this case we will say that both clusters are equivalent.
In order to have a well defined weighted cluster associated to every
cluster scheme, one considers the notion of \emph{consistent}
clusters, which we define next.
The \emph{proximity inequality} at $p_i$ is
$$ m_i \geq \sum_{p_j \text{ prox. to } p_i} m_j.$$
A weighted cluster $(K,\m)$ is consistent if and only if
it satisfies the proximity inequalities at
all its points.
Given a (unibranched) cluster scheme $Z$, there
is a unique consistent weighted cluster $(K,\m)$
such that $Z=Z_{K,\m}$ and $m_{i}>0$ for all $i$.
Furthermore, for a weighted cluster $(K,\m)$ non
necessarily consistent,
$$
\len Z_{K,\m}=\dim {\O \over {H_{K,\m}}} \leq
   \sum_{i=1}^r {{m_i\,(m_i+1)}\over 2} \, ,
$$
with the equality holding if $(K,\m)$ is consistent.

Given an arbitrary weighted cluster $(K,\m)$
there is a procedure called \emph{unloading}
(see \cite[4]{cll}, \cite[IV.II]{enr}, or \cite{cma})
which gives a new system of multiplicities
$\delta(\m)=\delta_K(\m)$
such that $(K,\delta(\m))$ is consistent and 
equivalent to $(K,\m)$. In each step of the
procedure, one \emph{unloads} some amount of multiplicity
on a point $p_i$ whose proximity inequality is not satisfied,
from the ponts proximate to it. This means
that there is an integer $n>0$ such that, increasing the
multiplicity of $p_i$ by $n$ and decreasing the
multiplicity of every point proximate to $p_i$
by $n$, the resulting weighted cluster
is equivalent to $(K,\m)$ and satisfies the
proximity inequality at $p_i$. In other words,
if $\tilde E_i \subset S_K$ is the strict transform
of the exceptional divisor of blowing-up $p_i$,
$D=-m_1 E_1 - m_2 E_2 - \cdots - m_r E_r$ and
and $\tilde E_i \cdot D <0$ then one chooses $n$ as
the minimal integer with $\tilde E_i \cdot (D-n \tilde E_i) \geq 0$
and replaces $D$ by $D-n \tilde E_i$.
A finite number of unloading steps lead to the
desired equivalent consistent cluster $(K,\delta(\m))$.

Let $C \subset S$ be a curve (more generally, a divisor).
For any proper or infinitely near point $q$ of $S$
we write $e_q(C)$ the multiplicity of (the strict
transform of) $C$ at $q$. Let $(K,\m)$ be a weighted cluster
of $r$ points
and $i \leq r$. Let $S_i$ be the surface
obtained by blowing up the first $i$ points of $K$, 
$\bar C \subset S_i$ the pullback of $C$ and for
$1 \leq j \leq i$, $E_j \subset S_i$ the pullback
of the exceptional divisor of blowing up $p_j(K)$.
We define the \emph{virtual transform} of $C$ in $S_i$
relative to the system $\m$ as the divisor
$$ \tilde C = \bar C -m_1 E_1 - m_2 E_2 - \cdots - m_i E_i \ .$$
A curve $C$ contains the cluster
scheme $Z_{K,\m}$ if and only if its virtual transform
in $S_{i-1}$ has multiplicity at least $m_i$ at $p_i$
for all $i$. Then we say that $C$ \emph{goes through}
$(K,\m)$. If $e_{p_i}=m_i$ for all $i$ then the virtual
transform coincides with the strict transform; if furthermore
$\bar C$ is a divisor with normal crossings
then we say that $C$ goes \emph{sharply} through $(K,\m)$.
In this case, $\pi_K$ is an embedded resolution of singularities
for $C$, and $(K,\m)$ determines the equisingularity class of $C$.
There are curves going sharply through $(K,\m)$ if
and only if $(K,\m)$ is consistent, and in this case they
are all equisingular.

We define next varieties $Y_i \subset X_i$
and smooth surjective morphisms
$\psi_i: X_{i} \rightarrow Y_{i-1}$
of relative dimension 2, as follows:
Let $Y_{-1}=\Spec k$, $X_0= S$ and $Y_0= \{ p \} $, and
for $i>0$, let
$$ X_i @>{b_i}>> Y_{i-1} \times_{Y_{i-2}} X_{i-1} $$
be the blowing-up along the diagonal
$  \Delta (Y_{i-1}) \subset Y_{i-1} \times_{Y_{i-2}} X_{i-1}$,
let $Y_i \subset X_i$ be the exceptional divisor, and
$\psi_i$ the composition
$\pi_{Y_{i-1}} \circ b_i :X_i \rightarrow Y_{i-1}$. We
define also morphisms $\pi_i=\pi_{X_{i-1}} \circ b_i$ and
$\pi_{i,j}=\pi_{j+1} \circ \dots \circ \pi_{i-1} \circ \pi_i:
X_i \rightarrow X_j$.

Remark that by construction the variety $Y_i$ is irreducible
and smooth, and $\psi_i|_{Y_i} :Y_i \rightarrow Y_{i-1}$
is a $\P^1$--bundle, therefore $Y_i$ is also
projective and rational for all $i$. The varieties $X_i$ are
irreducible, smooth, projective and rational if $S$ is.
$Y_{r-1}$ can be identified with the set of all
unibranched clusters of $r$ points of $S$ with origin
at $p$ (cf. \cite{meu}, \cite{hasun}, \cite{ran})
in such a way that
for all $K=(p_1, p_2, \dots, p_r) \in Y_{r-1}$,
denoting by $j_K: \{K \} \rightarrow Y_{r-1}$ the inclusion
 one has a pullback diagram
$$
\diagramcompileto{pullback}
   S_K \rrto^{i_K} \dto & & X_r \dto^{\psi_r} \\
   \Spec k  \rdouble & \{K\} \rto^{j_{K}} & Y_{r-1}  
\enddiagram
$$
where $i_K$ is a closed immersion, and a commutative square
$$
\diagramcompileto{explota}
   S_K \rrto^{i_K} \dto_{\sigma_r} & & X_r \dto^{\pi_r} \\
   S_{\breve K} \rrto^{i_{\breve K}} &&  X_{r-1} 
\enddiagram
$$
where
$\breve K=(p_1, p_2, \dots, p_{r-1}) =
 \psi_{r-1} (K) \in Y_{r-2}$ is the cluster
 of the first $r-1$ points
of $K$, and $\sigma_r$ is the blowing-up of $p_r(K)$.
Remark that $S_{\breve K} = S_{r-1}(K)$ and
$i_{\breve K} (p_r(K))=K$.

For every pair of integers $s, t$ such that
$r \geq s > t \geq 1 $
the subset of $Y_{r-1}$ containing exactly the
clusters $K$ with $p_s(K)$  proximate to $p_t(K)$
is an irreducible
closed subvariety $P_{s,t}$ (cf. \cite{meu}).
There are also open and dense subsets
$U_{s,t} \subset P_{s,t}$ which contain
all unibranched clusters with
$p_i(K)$ proximate to $p_{i-1}(K)$ for all $i>1$ and to $p_t(K)$
if $s \geq i > t$,
and no other proximity relations. We will write
$P_s=P_{s,1}$ and $U_s=U_{s,1}$.
We say that a system of multiplicities $\m$
is consistent in a subset $W \subset Y_r$ when the
weighted cluster $(K,\m)$ is consistent for all $K \in W$.

\begin{figure}
  \begin{center}
    \mbox{\includegraphics[2256,2480][2459,2666]{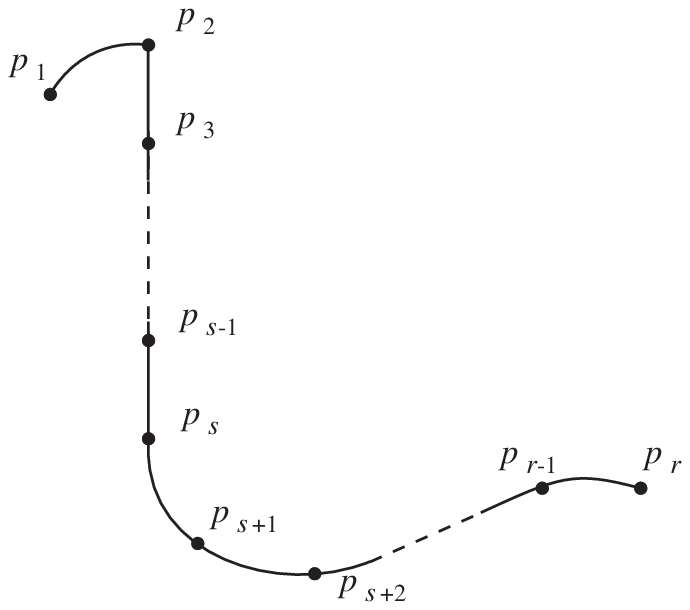}}
    \caption{Enriques diagram of a cluster in $U_s$}
  \end{center}
\end{figure}

\begin{Lem}
\label{filtra}
Fix a cluster $K=(p_1, p_2, \ldots, p_r) \in U_{s}$.
Let $i,j$ be positive integers such that $r=i+j+1$ and consider
the systems of multiplicities
\begin{align*}
\m_{-}&=(m_1,2^i,1^j)\\
\m_{+}&=(m_1,2^{i+1},1^{j-1}) \, .
\end{align*}
Let $I \subset \O$ be an ideal such that
$$ H_{K,\m_{+}} \varsubsetneq I \varsubsetneq H_{K,\m_{-}} \, .$$
Then there is a point $q \in E_r$ 
such that $I=H_{K_q,\m_{0}}$ with
$K_q=(p_1, p_2, \ldots, p_r, q)$ and
$\m_{0}=(m_1,2^i,1^{j+1})$ .
\end{Lem}

\begin{proof}
For any $q \in E_r$ let $q'$ be the unique point that is proximate
to $q$ and to $p_r$, and consider $K_q'=(p_1, p_2, \ldots, p_r, q, q')$
and $\m'=(m_1,2^i,1^{j+2})$. By performing successive unloading
steps on $p_r, p_{r-1}, \ldots, p_{i+2}$ we see that $(K_q',\m')$
is equivalent to $(K,\m_{+})$ and since $(K_q',\m')$ results
from $(K,\m_{-})$ by adding two simple points,
$$
\dim {{H_{K,\m_{-}}}\over{H_{K,\m_{+}}}} \leq 2 \, .
$$
The hypothesis on $I$ implies that the inequality is
in fact an equality and
$$
\dim {{I}\over{H_{K,\m_{+}}}} = \dim {{H_{K,\m_{-}}}\over{I}} =1 \, ,
$$
so $I=H_{K,\m_{+}}+(f)$ for some $f \in H_{K,\m_{-}} \setminus H_{K,\m_{+}}$.
$f$ is the equation of a germ of curve $\xi$ which goes through $(K,\m_{-})$
but not through $(K,\m_{+})$, therefore its virtual
transform at $p_r$ relative to the system of multiplicities
$\m_{-}$ is smooth, so it has a unique point
$q \in E_r$. It is clear that
$I=H_{K,\m_{+}}+(f) \subseteq H_{K_q,\m_{0}}$ and 
$\dim \left( {{H_{K_q,\m_{0}}} / {H_{K,\m_{+}}}} \right) \leq 1 $,
so $I=H_{K_q,\m_{0}}$.
\end{proof}

\begin{Pro}
\label{conductor}
Let $(K,\m)$ be a weighted cluster supported at $p$ and $f \in \O$ the
equation of a germ of curve $\xi$ through $p$. Write
$e_i=e_{p_i}(\xi)$. Then the conductor $(H_{K,\m}:f)$ is $H_{K,\m'}$
with $$m'_i=m_i - e_i \, .$$
\end{Pro}
\begin{proof}
The total transform of $\xi$ in $S_K$ is
$$ \tilde \xi + e_1 E_1 + e_2 E_2 + \cdots + e_r E_r $$
where $\tilde \xi$ is the strict transform. Let $g \in \O$ be
the equation of a second germ $\eta$. If $g \in H_{K,\m'}$ the
total transform of $\xi + \eta$ is
$$ \tilde \xi+ C + ((m_1-e_1)+e_1) E_1 +
      ((m_2-e_2)+e_2) E_2 + \cdots + ((m_r-e_r)+e_r) E_r $$
with $C$ effective, so clearly $fg \in H_{K,\m}$. Conversely,
if $fg \in H_{K,\m}$, then the total transform of $\xi + \eta$
is
$$ \tilde \xi+ C + m_1 E_1 + m_2 E_2 + \cdots + m_r E_r $$
with $C$ effective, so the total transform of $\eta$ must be
$$ C + (m_1-e_1) E_1 +
      (m_2-e_2) E_2 + \cdots + (m_r-e_r) E_r $$
and $g \in H_{K,\m'}$.
\end{proof}

If $Z \subset S$ is a zero--dimensional scheme defined
by the ideal sheaf $\I_{Z/S}$ and
$C \subset S$ is a curve, then there is an exact sequence
$$
0 \longrightarrow \I_{Z'/S} (-C)
  \longrightarrow \I_{Z/S}  \longrightarrow
  \I_{(Z \cap C) /C} \longrightarrow 0
$$
where $\I_{(Z \cap C) /C} = \I_{Z/S} \otimes \O_C$
and $\I_{Z'/S}$ defines a zero--dimensional scheme
$Z' \subset S$. In the context of the Horace method,
it is usually called the \emph{residual exact sequence}
and $Z'$ is the \emph{residual scheme} of $Z$ with respect
to $C$ (cf. \cite{hir}).

\begin{Cor}
\label{residual}
If $Z_{K,\m} \subset S$ is a cluster scheme and $C$ is
a curve with $e_{p_i}(C)=e_i$  then the residual scheme
$Z'$ of $Z$ with respect to $C$ is the cluster scheme
$Z_{K,\m'}$ with $\m'$ as in proposition \ref{conductor}.
\end{Cor}

\begin{Pro}
\label{isomorf}
Given a cluster $K=(p_1, p_2, \ldots, p_r) \in U_{s}$ and a point
$q \in E_r$, consider $K_q=(p_1, p_2, \ldots, p_r, q)$. Suppose
that $r \geq (s-1)(s-2)$. Then for any $q, q' \in E_r$, each of
them proximate only to $p_r$, there is an automorphism
$\varphi^{*} :\O \rightarrow \O$ such that, for any
system of multiplicities $\m$,
$$ \varphi^{*}(H_{K_q,\m})=H_{K_{q'},\m} \, .$$
In particular, $Z_{K_q,\m} \cong Z_{K_{q'},\m}$.
\end{Pro}
\begin{proof}
Let $\m_0=(s-1,1^r)$.
It will be enough to see that there are open neighbourhoods
$V_i \subset S_{i-1}$ of $p_i$ and $V_{r+1} \subset S_K$ containing
both $q$ and $q'$, and isomorphisms $\varphi_i: V_i \rightarrow V_i$
commuting with the blowing-ups, such that $\varphi_i(p_i)=p_i$
and $ \varphi_{r+1} (q)= q'$. This is equivalent to prove that
there are isomorphic unibranched germs of curve $\xi$ and $\xi'$
going sharphly through $(K_q,\m_0)$ and $(K_{q'},\m_0)$ respectively,
because then there is a neighbourhood $V$ of $p_1$ where both
$\xi$ and $\xi'$ have a representative and an automorphism
$\varphi$ of $V$ sending one to the other and which therefore
lifts to the desired $\varphi_i$.
Let $\xi: f=0$ and $\zeta: g=0$
be arbitrary unibranched germs of curve
going sharphly through $(K_q,\m_0)$ and $(K_{q'},\m_0)$.
They are equisingular with a single characteristic
exponent $s/(s-1)$ and their intersection multiplicity
at $p$ is $\nu\geq (s-1)^2 + (s-1)(s-2) -1$, so 
there is an automorphism of the completion
$$ \psi^{*}: \hat \O \rightarrow \hat \O $$
with $\psi^{*}(f)=g$ (cf. \cite{zpm}).
Let $x,y \in \O$ be a system of parameters;
we have $\hat \O \cong k[[x,y]]$, and $\psi^{*}$
can be described by
\begin{align*}
  \psi^{*}(x)&= \sum_{i,j>0} a_{ij}x^i y^j \\
  \psi^{*}(y)&= \sum_{i,j>0} b_{ij}x^i y^j 
\end{align*}
Choose $e \in \Z_{>1}$ such that $\M^{e} \subset H_{K_{q'},\m_0}$.
Then
\begin{align*}
  \varphi^{*}(x)&= \sum_{i+j\leq e} a_{ij}x^i y^j \\
  \varphi^{*}(y)&= \sum_{i+j\leq e} b_{ij}x^i y^j 
\end{align*}
define an automorphism of $\O$ such that
$\varphi^{*}(f)=0$ is the equation of an irreducible
germ $\xi'$ going sharply through $(K_{q'},\m_0)$.
Indeed, $f \in \M=(x,y)$ so
$\varphi^{*}(f)-\psi^{*}(f) \in \M^e \subset H_{H_{q'},\m_0}$
and we know also that $\psi^{*}(f) \in H_{H_{q'},\m_0}$.
Therefore $\varphi^{*}(f)=0$ is the equation of a
germ $\xi'$ going through $(K_{q'},\m_0)$. As
$\varphi^{*}$ is an automorphism, $\xi'$ is irreducible
and goes sharply through $(K_{q'},\m_0)$.
\end{proof}

\section{Specializing unibranched cluster schemes}
\label{especialitzar}

The subscheme of the Hilbert scheme which parametrizes
cluster schemes of a given type, and its relation to
the corresponding variety of clusters, has been studied
by many authors (\cite{nv}, \cite{pax}, \cite{ran}, \cite{gls2},
among others). However, they usually assume that the
proximity relations between points of the clusters
remain constant. Little seems to be known about
the relative position of these subschemes. Evain
computed in \cite{ev1} and \cite{ev3} several collisions of points,
including all cases with 3 points; this is
equivalent to the determination of the closure of
the corresponding subschemes of the Hilbert scheme,
that is, the specializations of cluster schemes when
new proximity relations arise. 
In this section we show some
flat families of unibranched cluster schemes in
which proximities vary; we introduce them in order
to prove theorems \ref{rang} and \ref{existencia}
but they can have also some interest on their own.

Let $x, y$ be local parameters for $\O$.
The inclusion $\O \subset \hat \O \cong k[[x,y]]$ allows
to write any $f \in \O$ as a formal power series
$$ f=\sum_{i,j \geq 0} a_{ij}x^i y^j $$
in a unique way. Thus any polynomial
$\Pi (X_{ij}) \in R=k[X_{ij}]_{i,j \geq 0}$
determines a function $\O \rightarrow k$, by
evaluating at $a_{ij}$: $\Pi(f)=\Pi(a_{ij})_{i,j \geq 0}$,
and to every $f \in \O$ corresponds a maximal ideal
$$ \M_f = \{ \Pi \in R \ | \ \Pi(f)=0 \} \subset R$$
with quotient field equal to $k$. So we have a mapping
$\O \rightarrow \MaxSpec R \subset \Spec R$, and it is
easy to see that it is injective. For any variety $Y$
we have therefore $\O \times Y \hookrightarrow \Spec R \times Y$,
and we will take the Zariski topology on $\O \times Y$
induced by the one on $\Spec R \times Y$.
%
%
It is easy to see that this Zariski topology does not
depend on the local parameters chosen for $\O$.
We will prove in \ref{tancat} below that for any system
of multiplicities $\m$ the set
$$H_\m:= \{ (f,K) \in \O \times Y_{r-1} \ | \  f\in H_{K,\m} \} \, .$$
is Zariski--closed in $\O \times Y_{r-1} \ .$

Fix for the rest of the section a system of multiplicities
$\m=(m_1, m_2, \dots, m_r)$ which we will apply to any cluster
of $r$ points; write $\breve \m=(m_1, m_2, \dots, m_{r-1})$. 
For a cluster $K=(p_1, p_2, \dots, p_r)$, we write also
$\breve K=(p_1, p_2, \dots, p_{r-1})$.

\begin{Lem}
\label{polinomis}
For every $K_0 \in Y_{r-1} \subset X_{r-1}$ there are an open
neighbourhood $V \subset X_{r-1}$ of $K$, an isomorphism
of $k$--varieties
$$U:= V \cap Y_{r-1} \overset{\lambda}\longrightarrow k^{r-1} \ ,$$
two functions
$$x, y \in \smash{\Gamma \left( U \times_{Y_{r-2}} V, \
   \O_{Y_{r-1} \times_{Y_{r-2}} X_{r-1}} \right) }$$
generating the
ideal of $\Delta(U)$ and polynomials
$A_{ij}^r \in R[u_1, u_2, \dots, u_{r-1}]$
for $i, j \geq 0$, such that for every $K \in U$ and
$f \in H_{K,\m}$ the formal power series
$$ \sum_{i,j \geq 0} A_{ij}^r(f,\lambda(K))
   (x \circ i_{\breve K})^i (y \circ i_{\breve K})^j \in \O_{S_{r-1},p_r}$$
and is a local equation for the virtual transform of $f=0$
at the last point $p_r$ of $K$.
\end{Lem}

Remark that the conditions on $x, y$ and $U$ imply that
for any cluster $K \in U$ the functions
$x, y, \lambda_1 \circ \pi_U - \lambda_1(K), \lambda_2 \circ
\pi_U - \lambda_2(K), \dots, \lambda_{r-1} \circ \pi_U - \lambda_{r-1}(K)$
are a system of parameters for the local ring
of $(K, K)$ in $U \times_{Y_{r-2}} V$,
the surface
$$\{K \} \times_{Y_{r-2}} X_{r-1} =
  \{K \} \times i_{\breve K}(S_{\breve K}) \cong S_{r-1}(K)$$
is locally defined by the equations $\lambda_i \circ \pi_U = \lambda_i(K)$,
and $x \circ i_{\breve K}$, $y\circ i_{\breve K}$ are local
parameters for $\O_{S_{r-1},p_r}$ (recall that
$S_{\breve K}=S_{r-1}(K)$ and $i_{\breve K} (p_r)= K$).

\begin{proof}
We proceed by induction on $r$. For $r=1$,
choose a system of parameters $x,y \in \O$ and a neighbourhood
$V$ of $p$ where $x$ and $y$ have regular representatives
(which abusing notation we call $x, y$ also) and
$V(x,y) \cap V=\{p\}$. Then the claim is clear,
because $Y_0=\{p\}$, and the virtual transform
of $f$ at $p$ is $f$ itself, so we can take $A_{ij}^1=X_{ij}$.
Suppose now $r>1$ and apply the induction hypothesis to
$\smash{(\breve K, \breve \m)}$. We obtain
the existence of $\breve V \subset X_{r-2}$,
$\breve U = \breve V \cap Y_{r-2} \cong k^{r-2}$,
$\breve \lambda$, $\breve x, \breve y$, and the polynomials $A_{ij}^{r-1}$
as in the claim.
By making a linear substitution in $\breve x, \breve y$ we may assume
that the last point $p_r$ of $K_0$ lies in the direction
of $\breve y \circ i_{\breve K_0}=0$.
$\smash{b_{r-1}^{-1}(\breve U \times_{Y_{r-3}} \breve V)}$ is open in
$X_{r-1}$, and is the blowing-up of $\smash{\Delta(\breve U)}$ in
$\smash{\breve U} \times_{Y_{r-3}} \breve V$, which can
be described as the subvariety of
$$\breve U \times_{Y_{r-3}} \breve V \times \P^1$$
given by the equation $u \breve x-v \breve y=0$, where $(u:v)$ are
projective coordinates of $\P^1$. The exceptional divisor
$Y_{r-1} \cap \smash{b_{r-1}^{-1}(\breve U \times_{Y_{r-3}} \breve V)}$
has equations
$\breve x= \breve y=0$. We define now $V$ to be the open subset
determined by $v \neq 0$, and the isomorphism
$$\lambda: U \longrightarrow k^{r-1}=k^{r-2} \times k$$
as $\lambda=(\breve \lambda \circ \psi_{r-1}) \times (u/v)$.
Let $\pi_V, \pi_U$ be the projections of
$V \times_{\breve U} U$ on its two factors. Then
it is easy to see that
\begin{align*}
x & = \breve x \circ \pi_V \\
y & ={u\over v} \circ \pi_V - {u \over v} \circ \pi_U =
   {\breve y \over \breve x} \circ \pi_V - \lambda_r \circ \pi_U
\end{align*}
generate the ideal of $\Delta(U)$.

Because of the induction hypothesis we know that for any
$K \in U$ and $f \in H_{K,\m}$ the virtual transform of
$f$ at $p_{r-1}(K)$ relative to $\m$ is
$$ \breve f = \sum_{i,j \geq 0} A_{ij}^{r-1}
   \left( f,\lambda \circ \psi_{r-1} (K) \right)
   \left( \breve x \circ i_{\breve {\breve K}} \right)^i
   \left( \breve y \circ i_{\breve {\breve K}} \right)^j .$$
This virtual transform must have
multiplicity at least $m_{r-1}$ at $p_{r-1}(K)$, therefore
$$ A_{ij}^{r-1}(f,\lambda \circ \psi_{r-1}(K)) =0 \
   \forall i,j, \ i+j < m_{r-1} \ ,$$
and the virtual transform transform of $\breve f$ in
$S_{r-1}(K)$ is given locally by
$$ \sum_{i+j \geq m_{r-1}} A_{ij}^{r-1}
   (f,\lambda \circ \psi_{r-1} (K))
   \left( \breve x \circ i_{\breve{\breve K}} \right) ^{i+j-m_{r-1}}
   \left( {\breve y \over \breve x} \circ i_{\breve {\breve K}}
   \right) ^j =0,$$
which because of the commutativity
$i_{\breve{\breve K}} \circ \sigma_{r-1} = \pi_{r-1} \circ i_{\breve K}$
can be written in terms of the local parameters
$x \circ i_{\breve K}, y \circ i_{\breve K}$:
\begin{gather*}
 \sum_{i+j \geq m_{r-1}} A_{ij}^{r-1}
   (f,\lambda \circ \psi_{r-1} (K))
   \left( x \circ i_{\breve K}\right) ^{i+j-m_{r-1}}
   \left( y \circ i_{\breve K} + \lambda_r(K) \right) ^j = \\
= \sum_{i,j \geq m_{r-1}} \sum_{\ell=0}^j A_{ij}^{r-1}
   (f,\lambda \circ \psi_{r-1} (K))\binom{j}{n} \lambda_r(K)^{j-\ell}
   \left( x \circ i_{\breve K} \right) ^{i+j-m_{r-1}}
   \left( y \circ i_{\breve K} \right) ^{\ell} \ .
\end{gather*}
This allows us to define the polynomials
$$ A_{k\ell}^r (f,u_1, \dots, u_{r-1})=
 \sum_{i+j=k+m_{r-1}} A_{ij}^{r-1} (f, u_1, \dots, u_{r-2})
   \binom{j}{n} u_{r-1}^{j-\ell} \ ,$$
after which the claim is satisfied.
\end{proof}

\begin{Pro}
\label{tancat}
$H_\m$ is Zariski--closed in $\O \times Y_{r-1}$.
\end{Pro}

\begin{proof}
By induction on $r$. The case $r=1$ is clear, because
then $H_{K,\m}=\M^{m_1}$ is the closed subset of
$\O=\O \times Y_0$ determined by the ideal
$I=(X_{ij})_{i+j<m_1}$.
Suppose now $r>1$ and the claim is true for
$H_{\breve \m} \subset \O \times Y_{r-2}$.
Let $K_0 \in Y_{r-1}$ be a cluster and
$U \subset Y_{r-1}$ the open neighbourhood given by lemma
\ref{polinomis}. It will be enough to see that
$H_{U,\m}:=H_\m \cap \pi_{Y_{r-1}}^{-1}(U)$ is closed
in $\O \times U$, because $Y_{r-1}$ can be covered by
a finite number of such open neighbourhoods.
Define
$$
H'_{U, \m} = \{ (f,K) \in \O \times U \ | 
 (f, p_r(K)) \in H_{\breve \m} \} \ .
$$
Because of the induction hypothesis
and the fact that
$$ id_\O \times p_r: \O \times Y_{r-1} \longrightarrow
   \O \times Y_{r} $$
is continuous, $H'_{U,\m}$ is closed in $\O \times U$.
Moreover $(f,K) \in H_{U,\m}$ if and only if $(f,K) \in H'_{U,\m}$
and the virtual transform of $f=0$ at $p_r(K)$ has multiplicity
at least $m_r$, that is, by lemma \ref{polinomis},
$$ A_{ij}^{r}(f,\lambda (K)) =0 \
   \forall i,j, \ i+j < m_{r} \ .$$
These equations define a closed subset
$V(A_{ij}^r)_{i+j<m_r} \subset \O \times U$,
therefore $H_{U,\m}= H'_{U,\m} \cap V(A_{ij}^r)_{i+j<m_r}$
is closed in $\O \times U$.
\end{proof}

From this we obtain a number of corollaries.
Let $N$ be a positive integer such that
$\dim ( \O / H_{K,\m} ) \leq N$ for all $K\in Y_{r-1}$
(for example, $N=\sum_{i=1}^r {m_i\,(m_i+1)}/ 2$).
Then $H_{K,\m} \supset \mathfrak m^{N}$ for all $K$, and
we can define
$$ \bar{H}_\m:=\{ (\bar f,K) \in {\O \over {\M^N}}\times Y_{r-1} \
| \  f\in H_{K,\m} \} \, .$$
\begin{Cor}
$\bar{H}_\m$ is Zariski--closed in $\left( \O / \M^N \right) \times Y_{r-1}$.
\end{Cor}
\begin{proof}
Let $K_0 \in Y_{r-1}$ be a cluster and
$U \subset Y_{r-1}$ the open neighbourhood given by lemma
\ref{polinomis}. As in lemma \ref{tancat}, it will be enough
to see that $\bar H_{U,\m}:=\bar H_\m \cap \pi_{Y_{r-1}}^{-1}(U)$
is closed in $( \O / \M^N) \times U$. We also know from
the proof of lemma \ref{tancat} that $H_{U,\m}$ is defined
by a finite number of polynomials $A_{ij}^k$, $k=1, \dots, r$.
Now for $\alpha + \beta \geq N$,
$x^{\alpha}y^{\beta}\in H_{K,\m}$ for all $K$,
so the polynomial 
$A_{ij}^k(x^{\alpha}y^{\beta},u_1, \dots, u_{k-1})
  \in k[u_1, \dots, u_{k-1}]$
must be identically
zero. This implies in fact that
$$A_{ij}^k \in k[X_{\alpha \beta}]_{\alpha+\beta <N}
   [u_1, u_2, \dots, u_{k-1}]
 \subset R[u_1, u_2, \dots, u_{k-1}] \ .$$
Now $k[X_{\alpha \beta}]_{\alpha+\beta <N}$ is the
affine coordinate ring of $\O / \M^N$, which
is a $k$--vector space admiting the basis
$\{\bar x^\alpha \bar y^\beta \}_{\alpha+\beta <N}$.
Therefore the $A_{ij}^k$ define a Zariski--closed subset
in $\left( \O / \M^N \right) \times U$, and it is immediate
to see that this is in fact $\bar{H}_{U,\m}$.
\end{proof}

\begin{Cor}
\label{semic}
The function
\begin{align*}
Y_{r-1} &\longrightarrow {\Bbb Z}_{\geq 0} \\
 K &\longmapsto \dim {{\O \over H_{K,\m}}} = \len Z_{K,\m}
\end{align*}
is lower semicontinuous.
\end{Cor}
\begin{proof}
Choose $N$ such that 
$\dim ( \O / H_{K,\m} ) \leq N$ for all $K\in Y_{r-1}$.
Then
$$\dim {{\O \over H_{K,\m}}} =
      \dim {{\O / \M^N} \over {H_{K,\m}/ \M^N}} =
      {{N(N+1)} \over 2} - \dim \left(\bar H_{\m}
      \times_{Y_{r-1}} \{K\} \right)$$
and the claim follows because the dimension of the fibers
of the (non flat) family $\bar H_{\m} \rightarrow {Y_{r-1}}$
is upper semicontinuous.
\end{proof}

Remark that in general there is no ideal sheaf $\I$
on $S \times Y_{r-1}$ with $\I \otimes k(K)= H_{K,\m}$
for all $K$; in other words, all $Z_{K,\m}$ do not form
a family. Indeed, the length of the members of a family
of zero--dimensional schemes is upper semicontinuous,
against corollary \ref{semic}. As the simple example
$\m=(2,2,2)$ shows, this length is not always constant,
and in fact the systems of multiplicities for which it is
constant are quite exceptional. However, we will see that
restricting the set of clusters under consideration so that
the length of the corresponding schemes remains constant,
it is possible to construct flat families of cluster schemes.

Given a subvariety $W \subset Y_{r-1}$, let
$W_{\m} \subset W$ be the open subset where
$\dim \left(\O / H_{K,\m} \right)$ is maximal, and
denote this dimension by $N=N(\m,W)$.
The restriction
$$\bar H_\m |_{W_{\m}} \subset {\O \over {\M^N}}\times W_\m$$
is thus a family of vector subspaces of $\O/{\mathfrak m}^{N}$
of codimension $N$, so it defines a morphism
$W_\m \rightarrow \G$ to the
Grassmannian ${\G}={\G_{N}}\left( \O / \M^{N} \right)$.
In $\G$, the codimension $N$ vector spaces which are
ideals of $\left( \O / \M^{N} \right)$ constitute a closed
subscheme, which is identified to the 
Hilbert scheme $\Hilb^{N}_p S$ of length $N$ subschemes of
$S$ supported at $p$ (cf. \cite{fog}, \cite{iar}).
As $H_{K,\m}$ is by definition an ideal for every $K$,
we finally obtain a morphism
\begin{align*}
 W_{\m} &\overset{\varphi_\m}\longrightarrow \Hilb^{N} S \\
   K &\longmapsto Z_{K,\m}
\end{align*}
That is, the schemes $Z_{K,\m}$ form a flat family
with parameter space $W_{\m}$.
In the particular case that $W=C$
is a smooth curve, the morphism $\varphi_\m$
can be extended to all of $C$, by \cite[III, 9.8]{hag},
but then it is not true that $\varphi_\m (K)=Z_{K,\m}$
for those $K \in C \setminus C_{\m}$. In fact, $\varphi_\m (K)$
depends not only on $K$ and $\m$ but also on $C$.
Nevertheless, we have the following

\begin{Cor}
\label{comparalimits}
If $C \subset Y_{r-1}$ is a smooth curve and
$K \in C$ then $Z_{K,\m} \subseteq \varphi_\m (K)$.
\end{Cor}
\begin{proof}
The morphism $\varphi_\m: C \rightarrow \Hilb^N_p S \subset \G$
defines a closed subset $\Phi_\m \subset (\O / \M^N) \times C$,
and we have to see that $\Phi_\m \subset \bar H_\m|_C$. As
$\bar H_\m$ is closed and
$$
\Phi_\m |_{C_\m}= \bar H_\m |_{C_\m} \subset {\O \over {\M^N}} \times C_\m
$$
it will be enough to see that $\Phi_\m$ is the closure
of $\Phi_\m|_{C_\m}$.

Take $(\bar f, K_0) \in \Phi_\m$, with $K_0\in C \setminus C_\m$.
The family determined by $\varphi_\m$ corresponds to an ideal
sheaf $\I \subset (\O / \M^N) \otimes \O_C$. By definition,
$\bar f \in \I \otimes_{\O_C} k(K_0)$, so there is an open neigbourhood
$U$ of $K_0$ in $C$ and a section
$$\bar f_U \in \Gamma
   \left( \Spec \frac{\O}{\M^N} \times U , \I \right)$$
with $\bar f= \bar f_U (K_0)$. As $U$ is open,
we may assume $U \setminus \{ K_0 \} \subset C_\m$, therefore
$\forall K \in U \setminus \{ K_0 \}, \ \bar f_U (K)
   \in \I \otimes k(K) = \bar H_{K,\m}$.
So we have a morphism
\begin{align*}
U & \longrightarrow \Phi_\m \\
K & \mapsto (\bar f_U (K), K)
\end{align*}
whose image contains $(\bar f, K_0)$ and has every other point
in $\Phi_\m|_{C_\m}$, therefore $(\bar f, K_0)$ lies in the
closure of $\Phi_\m|_{C_\m}$. As this holds for every
$(\bar f, K_0) \in \Phi_\m$ and $\Phi_\m$ is closed, the
claim is proved.
\end{proof}

Next lemma shows an example of a flat family of cluster
schemes which will be useful later on.

\begin{Lem}
\label{unpuntmes}
Let $(K,\m)$ be a consistent weighted cluster of $r$
points and $E=E_r \subset S_K$ the exceptional divisor
of blowing up the last point. Consider the system
of multiplicities $\m_1=(m_1, m_2, \dots, m_r, 1)$
and for every $q \in E$, $K_q=(p_1, p_2, \dots, p_r, q)$.
Then the schemes $Z_{K_q,\m_1}$ form a flat family.
\end{Lem}

\begin{proof}
It is clear that
$\len Z_{K_q,\m_1} \le Z_{K,\m} +1 \ \forall q$.
On the other hand, as $(K, \m)$ is consistent,
there are curves going sharply through $(K, \m)$
which miss $q$ (cf. \cite[4]{cll}), so
$\len Z_{K_q,\m_1} = Z_{K,\m} +1 \ \forall q$.
Abusing slightly notations, we will call
$E=i_K(E) \subset Y_r$; as the length of
$Z_{K_q,\m_1}$ does not depend on $q$ we have
$E_{\m_1}=E$ and the schemes $Z_{K_q,\m_1}$ form a flat family.

\end{proof}

Recall that $P_s \subset Y_{r-1}$ is the irreducible
subvariety containing the unibranched clusters $K$ with
$p_{s}(K)$ proximate to $p_1$, and $U_s \subset P_s$
is the open dense subset which contains all unibranched
clusters based at $p_1$ with $p_i(K)$ proximate to
$p_{i-1}(K)$ for all $i>1$ and to $p_1(K)$ if $s \geq i > 1$,
and no other proximity relations. Assume from now
that the system of multiplicities we are dealing with
is of the form $\m=(m,2^i,1^j)$. We will denote by
$ \V(s,m,i,j) = \overline{\varphi_{\m}(U_s)} $
the closure of the subscheme of the Hilbert scheme that
parametrizes the unibranched cluster schemes $Z_{K,(m,2^i,1^j)}$
with $K \in U_{s}$. Note that this does not depend on $r$,
as far as $r \geq i+j+1$. Note also that $U_s \subset (P_s)_\m$
because $U_s$ is open and dense in $P_s$ and
$\len Z_{K,\m}$ is constant on $U_s$. Therefore
$\varphi_\m ((P_s)_\m) \subset \overline{\varphi_\m (U_s)}$.


\begin{Pro}
\label{limit}
Suppose that $\m=(m,2^i,1^j)$ is consistent
in $U_{s}$ and $(P_s)_{\m}$ does not contain $U_{s+1}$.
Then
\begin{enumerate}
\item $m=2\,s -2$ and $i \geq s$,
   \label{primera}
\item if $i+j \geq s^2-3\,s+1$ then $\V(s,m,i,j)$ contains
   $\V(s+1,m+1,i-s+1,j+s-2)$,
\item $(m+1,2^{i-s+1}, 1^{j+s-2})$ is consistent in $U_{s+1}$ if and only if
$i \leq 2\,s-2$. If it is not, then it is equivalent to
$(m+2, 2^{i-2s+1}, 1^{j+2s-2})$.
\end{enumerate}
\end{Pro}

\begin{figure}
  \begin{center}
    \mbox{\includegraphics[2135,2486][2469,2674]{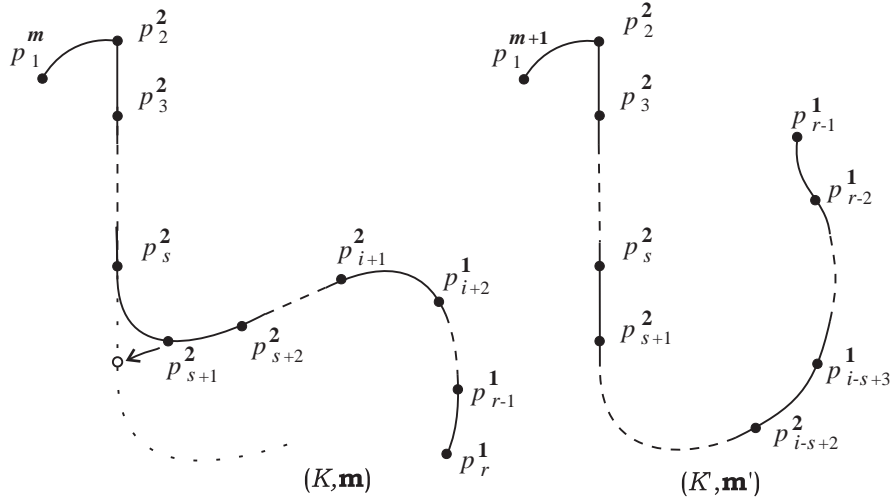}}
    \caption{Enriques diagrams illustrating proposition
      \ref{limit}. By moving the free points, $Z_{K,\m}$
      specializes to $Z_{K',\m'}.$}
  \end{center}
\end{figure}

\begin{proof}
\begin{enumerate}
\item $\m$ is consistent in $U_s$, but 
not consistent in $U_{s+1}$, otherwise
$U_{s+1} \subset (P_s)_{\m}$. Therefore the proximity
inequality at $p_1$ must be satisfied for $K \in U_s$
but not for $K \in U_{s+1}$. This means that
$$
\begin{cases}
2\,s > m \geq 2\,s -2 & \text{ if } i \geq s, \\
s+i > m \geq s+i-1 & \text{ if } i < s.
\end{cases}
$$
Therefore, if $i < s$ we must have $m=s+i-1$,
and for any cluster $K \in U_{s+1}$ we have
$\delta(\m)=(m+1,1^{2i+j-s})$. So
\begin{gather*}
\len Z_{K,\m}=\len Z_{K,\delta(\m)}= \\
={{(m+1)(m+2)} \over 2} +2\,i +j -s=
{{m\,(m+1)} \over 2} +3\,i +j
\end{gather*}
for all $K \in U_{s+1}$, so $U_{s+1} \subset (P_s)_\m$
against the hypotheses. We conclude that $i \geq s$, and it
only remains to be seen that the case $m=2\,s-1$ is not
possible. But in this case 
$\delta(\m)=(m+1,2^{i-s},1^{j+s})$ in $U_{s+1}$, so again
\begin{gather*}
\len Z_{K,\m}=\len Z_{K,\delta(\m)}= \\
={{(m+1)(m+2)} \over 2} +3\,(i-s) +j +s=
{{m\,(m+1)} \over 2} +3\,i +j
\end{gather*}
for all $K \in U_{s+1}$ and $U_{s+1} \subset (P_s)_\m$.

\item We are given two systems of multiplicities
\begin{align*}
\m&=(2\,s-2, 2^{i}, 1^{j}) \\
\m'&=(2\,s-1, 2^{i-s+1}, 1^{j+s-2})
\end{align*}
with $i>s$ and $i+j+1 \geq (s-1)(s-2)$. We have to see that for any
cluster $K_0 \in U_{s+1}$ there is a deformation of $Z_{K_0,\m'}$
whose general member is of the form $Z_{K,\m}$ with $K\in U_s$.
For that, let $C$ be a smooth irreducible curve $K_0 \in C \subset P_s$
with $C \cap U_s \neq \emptyset$ and consider the morphism
$$ C \overset{\varphi_\m}\longrightarrow \Hilb^N S \, .$$
Consider also an auxiliary system of multiplicities
and the associated morphism
\begin{gather*}
\m_{+}=(2s-1, 2^{i-s+1}, 1^{j+s-1}) \\
C \overset{\varphi_{\m_{+}}}\longrightarrow \Hilb^{N+1} S \, .
\end{gather*}
For $K \in C \cap U_s$, $Z_{K,\m} \varsubsetneq Z_{K,\m_{+}}$.
Therefore $\varphi_\m(K) \varsubsetneq \varphi_{\m_{+}}(K)$ for
all $K \in C$. Moreover, part \ref{primera} of the
proposition shows that $K_0 \in C_{\m_{+}}$, so
$\varphi_{\m_{+}}(K_0)=Z_{K_0,\m_{+}}$.
On the other hand, $(K_0,\m)$ does not satisfy the proximity
inequality at $p_1$, and unloading multiplicity on
this point gives
$$
\m_{-}=(2\,s-1, 2^{i-s}, 1^{j+s}) \,.
$$
Note that $\m_{-}$ may be non--consistent, but still
$Z_{K_0,\m}=Z_{K_0,\m_{-}}$ and
$\len Z_{K_0,\m}< \len \varphi_\m(K_0)$. All together, we have
$$
Z_{K_0,\m_{-}}=Z_{K_0,\m} \varsubsetneq \varphi_\m(K_0)
    \varsubsetneq \varphi_{\m_{+}}(K_0) = Z_{K_0,\m_{+}} 
$$
and we are in the conditions of lemma \ref{filtra}.
Therefore there is a point $q$ in the exceptional
divisor $E=E_r$ of blowing up $p_r$
such that $\varphi_{\m}(K_0)=Z_{K_q,\m_0}$, with
$$\m_0=(2\,s-1,2^{i-s},1^{j+s+1})\, .$$
Let $q_0 \in E$ be the only point proximate to $p_{r-1}$.
If $q=q_0$ then unloading gives
$Z_{K_q,\m_0}=Z_{K_0,\m'}$ so $\varphi_{\m}$ is the
family we are looking for.
If $q$ is not proximate to $p_{r-1}$ then we only know
that $Z_{K_q,\m_0} \in \V(s,m,i,j)$. 
As $i+j+1 \geq (s-1)(s-2)$, though, by lemma \ref{isomorf}
we can say that for any free point $q \in E$,
$Z_{K_{q},\m_0} \in \V(s,m,i,j)$. As the free
points are dense in $E$, this implies
$\varphi_{\m_0}(E) \subset \V(s,m,i,j)$. Now it
is enough to see that $q_0 \in E_{\m_0}$, because then
$Z_{K_0,\m'}=Z_{K_{q_0},\m_0}=\varphi_{\m_0}(q_0) \in \V(s,m,i,j)$.
If $(K_q,\m_0)$ is consistent for $q \in E$ free,
then $q_0 \in E_{\m_0}$ because of lemma \ref{unpuntmes};
if it is not, then the equivalent consistent system
obtained by unloading is
$\delta_{K_q}(\m_0)=(m+2, 2^{i-2s}, 1^{j+2s+1})$, whereas
for the cluster $K_{q_0}$ the equivalent consistent system
is $\delta_{K_{q_0}}(\m_0)=(m+2, 2^{i-2s+1}, 1^{j+2s-2})$.
This implies that $\len Z_{K_q,\m_0}= N(\m_0,E)$, so
again $q_0 \in E_{\m_0}$.
\item Follows from an easy unloading calculation.
\end{enumerate}
\end{proof}

\begin{Cor}
Suppose that $\m=(m,2^i,1^j)$ 
is consistent in $U_{s}$ 
and $(P_s)_{\m}$ does not contain $U_{s+1}$, and
suppose that $j\geq s^2-5\,s+2$.
For every $k \in \Z_{\geq 0}$ we define
\begin{align*}
m_k&=m+2\,k \\
i_k&=i-k\,(k+2\,s-2) \\
j_k&=j+k\,(k+2\,s-3) \, .
\end{align*}
If either $k=0$ or $i_{k-1} > 2\,(s+k-1)-2$ then
$\V(s+k,m_k,i_k,j_k) \subset \V(s,m,i,j)$.
\end{Cor}

\begin{proof}
We will proceed by induction on $k$. For $k=0$,
there is nothing to prove. For $k>0$, we have
either $k-1=0$ or
$$i_{k-2} > i_{k-1} > 2\,(s+k-1)-2 > 2\,(s+k-2) -2$$
so we can apply the induction hypothesis and
$$\V(s+k-1,m_{k-1},i_{k-1},j_{k-1}) \subset \V(s,m,i,j)\, .$$
On the other hand, a straightforward computation
shows that $j_{k-1}\geq (s+k-1)^2 -5\,(s+k-1)+2$,
therefore $i+j \geq (s+k-1)^2 -3\,(s+k-1)+1$.
In these contitions, proposition \ref{limit} tells us that
$$\V(s+k,m_{k},i_{k},j_{k}) \subset 
      \V(s+k-1,m_{k-1},i_{k-1},j_{k-1}) \, ,$$
thus finishing the proof.
\end{proof}

\section{Degree of singular plane curves}

In this section we work on a projective irreducible
smooth surface $S$, mainly $S=\P^2$. We are interested
in linear systems of curves which contain a zero
dimensional scheme composed of unibranched cluster schemes.
These linear systems can be specialized to linear systems of
curves through a unibranched cluster (supported at a single point)
using the technique developped in \cite{meu}. We will
not repeat the whole treatment here, but only show the
parameter space (a variety of clusters) to which we apply
the semicontinuity theorem.

For convenience, we fix a point $p$ in $S$.
Given a sequence of integers $\r=(r_1, r_2, \ldots, r_{\rho})$
with $r_1 < r_2 < \ldots < r_{\rho}$ we define
varieties $X_i(\r)$ and
$Y_i(\r)$ in an analogous way to
$X_i$ and $Y_i$. Let $Y_{-1}(\r)=\Spec k$, $X_0(\r)=S$,
$Y_0(\r)=\{p\}$, and for $i>0$ let
$$
X_i(\r) @>{b_i}>> Y_{i-1}(\r) \times_{Y_{i-2}(\r)} X_{i-1}(\r)
$$
be the blowing-up along $\Delta(Y_{i-1}(\r))$,
$$ Y_i(\r)=
\begin{cases}
\text{exc. divisor of the blowing-up} & \text{if } i \notin
   \{r_1, r_2, \ldots, r_{\rho}\} \\
X_i(\r) & \text{if } i \in
   \{r_1, r_2, \ldots, r_{\rho}\} 
\end{cases} $$
and $\psi_i(\r)=\pi_{Y_{i-1}(\r)} \circ b_i:
   X_i(\r) \rightarrow Y_{i-1}(\r)$
is smooth of relative dimension 2. We write also
$X(\r)=X_{r_{\rho}}(\r)$ and
$Y(\r)= Y_{r_{\rho}-1}(\r)$.
All these varieties are projective, irreducible
and smooth; $Y(\r)$ can be identified with
the set of all clusters of $r_{\rho}$ points of $S$, such that
$p_1=p$ and $p_{i+1}$ is proximate to $p_i$ for all
$i \notin \r$.
A general cluster $K \in Y(\r)$ is the
union of $\rho$ unibranched clusters of $r_1, r_2-r_1,
\ldots,$ and $r_{\rho}-r_{\rho -1}$ points.

\begin{Lem}
\label{ajuntar}
For any sequence of integers $\r=(r_1, r_2, \ldots, r_{\rho})$,
$X_i$ and $Y_i$ are closed subvarieties of $X_i(\r)$ and $Y_i(\r)$
respectively, and the morphism $\psi_i:X_i \rightarrow Y_{i-1}$
of section \ref{preliminars} is the restriction of 
$\psi_i(\r):X_i(\r) \rightarrow Y_{i-1}(\r)$ for all $i$.
\end{Lem}
\begin{proof}
By induction on $i$. The cases $i \leq 1$ are clear, so
assume $i>1$ and the claim to be true for $X_{i-1}$ and
$X_{i-2}$. This means that $X_{i-1} \subset X_{i-1}(\r)$,
$Y_{i-1} \subset Y_{i-1}(\r)$ are closed subvarieties
and the morphism $\psi_{i-1}:X_{i-1} \rightarrow Y_{i-2}$
is the restriction of
$\psi_{i-1}(\r):X_{i-1}(\r) \rightarrow Y_{i-2}(\r)$.
So we have also a closed subvariety
$$
Y_{i-1} \times_{Y_{i-2}} X_{i-1} \subset
Y_{i-1}(\r) \times_{Y_{i-2}(\r)} X_{i-1}(\r)
$$
and 
$$
\Delta(Y_{i-1}(\r)) \cap Y_{i-1} \times_{Y_{i-2}} X_{i-1}
= \Delta(Y_{i-1})
$$
so by the definitions, $X_i \subset X_{i-1}(\r)$ is the
strict transform of $Y_{i-1} \times_{Y_{i-2}} X_{i-1}$
under the blowing-up $b_i$ (cf. \cite[II.7.15]{hag}),
$Y_i$ is a subvariety of the exceptional divisor,
hence of $Y_{i-1}(\r)$, and $\psi_i:X_i \rightarrow Y_{i-1}$
is the restriction of $\psi_i:X_i(\r) \rightarrow Y_{i-1}(\r)$.
\end{proof}

Given a cluster $K \in Y(\r)$, a 
system of multiplicities $\m$ and a divisor $D$, we will
denote by $\L_D(K,\m)$ the linear system of effective divisors
linearly equivalent to $D$ which go through $(K,\m)$.
As seen in \cite{meu} the function
\begin{align*}
Y(\r) &\longrightarrow {\Bbb Z}_{\geq -1} \\
 K &\longmapsto \ell_D(K,\m) := \dim \L_D(K,\m)
\end{align*}
is upper semicontinuous. We consider the dimension of the
empty linear system to be $-1$. Because of lemma \ref{ajuntar},
 we can bound $\ell_D(K,\m)$ with $K$ general in
$Y(\r)$ by $\ell_D(K,\m)$ with $K$ in $Y_{r_{\rho}-1}$.
For any subset $W \subset Y(\r)$ we will also write
$\ell_D(W,\m)=\inf \{\ell_D(K,\m)\ | \ K\in W \}$.

In the case $S=\P^2$, we can choose $D=d\cdot L$, with
$L$ a line and $d \in \Z_{>0}$. We will then denote
$\L_d(K,\m)=\L_D(K,\m)$
the linear system of curves of degree $d$ going through
$(K,\m)$ and $\ell_d(K,\m)=\ell_D(K,\m)$ its dimension.
We say that a zero--dimensional scheme $Z$ has level $d$
when $\len Z =(d+1)(d+2)/2$. For a level $d$ scheme in the
plane $\P^2$ it is equivalent to have maximal rank or to
have maximal rank in degree $d$. Furthermore, an arbitrary
zero--dimensional subscheme $Z \subset \P^2$ has maximal
rank if and only if there exist maximal rank 
schemes $Z_d$ and $Z_{d+1}$ of level $d$ and $d+1$ respectively
such that $Z_d \subset Z \subset Z_{d+1}$ (see \cite[2.2.2]{hir}).

\begin{Lem}
\label{nivelld}
For every cluster scheme $Z_{K,\m}$ where $K \in U_s$ and
$\m=(m,2^i, 1^j)$ is consistent in $U_s$ there exist systems
of multiplicities $\m_{-}=(m,2^{i_{-}}, 1^{j_{-}})$ and
$\m_{+}=(m,2^{i_{+}}, 1^{j_{+}})$ consistent in $U_s$ such
that $Z_{K,\m_{-}} \subset Z_{K,\m} \subset Z_{K,\m_{+}}$,
$Z_{K,\m_{-}}$ has level $d$ and $Z_{K,\m_{+}}$ has level
$d+1$. Furthermore, if $4j \geq m^2 - 4m -6$ and
$3i+j \geq 2m+3$ then $\m_{+}$ and $\m_{-}$ can be
chosen such that $4j_{+} \geq m^2 - 4m -6$,
$4j_{-} \geq m^2 - 4m -6$ or $i_{-}=0$, and $d>m$.
\end{Lem}

\begin{proof}
Write
$$
{{m(m+1)} \over 2} + 3\, i + j = {{(d+1)(d+2)} \over 2}+ \epsilon
$$
with $0 \leq \epsilon \leq d+1$. It is enough to define
$\m_{+}=(m,2^i,1^{j+d+2-\epsilon})$ and
$$
\m_{-} = \begin{cases}
(m, 2^{i-\epsilon}, 1^{j+2\epsilon}) & \text{if } i \geq \epsilon \\
(m, 1^{j+3i -\epsilon}) & \text{if } i \leq \epsilon
\end{cases}
$$
\end{proof}

Before we prove our maximal rank theorem we need some lemmas
on level $d$ cluster schemes. Assume for a while that
$\m=(m,2^i,1^j)$ and $d \in \Z$ are such that
$$N(\m)={{m\,(m+1)}\over 2} + 3\,i +j = {{(d+1)\,(d+2)}\over 2}\,.$$

\begin{Lem}
\label{totuns}
If $i=0$ and $s \leq m+1$ then $\ell_d(U_s,\m)=-1$, with
$U_s \subset Y_{r_{\rho}-1}$ as defined in section \ref{especialitzar}.
\end{Lem}
\begin{proof}
It is clear that $d \geq m-1$; we will prove the claim by
induction on $d-m$. For $d=m-1$ the result is obvious.
For $d\geq m$ we have $j \geq m+1$. By semicontinuity
it is enough to see that there are no plane curves
of degree $d$ containing
$Z_{K,\m}$ with $K$ general in $U_{m+1}$. Now for
$K \in U_{m+1}$, unloading gives $Z_{K,\m}=Z_{K,(m+1,1^{j-m-1})}$,
and the result follows from the induction hypothesis.
\end{proof}

\begin{Lem}
\label{foradosos}
If $2\,i \leq m$, $2\,i +j>m$ and $\m$ is consistent in $U_s$
then $\ell_d(U_s,\m)=-1$.
\end{Lem}
\begin{proof}
By semicontinuity it is enough to see that there are
no plane curves of degree $d$ containing $Z_{K,\m}$ with $K$ general
in $U_{m+1-i}$. But for $K \in U_{m+1-i}$, unloading
gives $Z_{K,\m}=Z_{K,(m+1,1^{2i+j-m-1})}$,
and the result follows from lemma \ref{totuns}.
\end{proof}

\begin{Pro}
\label{horace}
If $d=m+1$, $m \geq 2i-2$ and $\m$ is consistent in $U_s$ then
$\ell_d(U_s,\m)=-1$.
\end{Pro}
\begin{proof}
The proof runs by induction on $i$, the case $i=0$ being clear
from lemma \ref{totuns}. Assume $i \geq 1$. For $K$ general in $U_s$,
consider the line $L$ passing through $p$ in the direction of $p_2$.
If $s=2$ we can assume ($K$ being general) that $p_3$ does not
lie on $L$; for $s>2$ this is automatic.
For every curve $C \in \L_d(K,\m)$
$$ L \cdot C \geq m+2 = d+1 > \deg C,$$
therefore $L$ is a component of $C$, so it is a fixed part of
$\L_d(K,\m)$. The residual linear system is
$$\L_d(K,\m) - L = \L_{d-1}(K,\m')$$
with $\m'=(m-1,1,2^{i-1},1^j)$ because of lemma \ref{residual}.
As this is true for general $K \in U_s$ we have now
$\ell_d(U_s,\m)=\ell_{d-1}(U_s,\m')$. But $\m'$ is not consistent
in $U_s$. The equivalent consistent system is
$$
\delta(\m')=
\begin{cases}
(m,1^{j+i}) & \text{ if } m= 2\,(s-1)=2\,(i-1) \\
(m-1,2^{i-1},1^{j+1}) & \text{ in any other case.}
\end{cases}
$$
Observe that $N(\delta(\m))=d(d+1)/2$. Finally
$\ell_d(U_s,\m) = \ell_{d-1}(U_s,\delta(\m'))=-1$ because of
lemma \ref{totuns} in the first case and because of
the induction hypothesis in the second.
\end{proof}

\begin{Cor}
\label{ipetit}
If $2i \leq m$, $d>m$ and $\m$ is consistent in $U_s$
then $\ell_d(U_s,\m)=-1$.
\end{Cor}
\begin{proof}
By lemma \ref{foradosos} we may assume $2i + j \leq m$.
But then
$$ N(\m)={{m(m+1)}\over 2} + 3\,i +j \leq
        {{m(m+1)}\over 2} + m +{m\over 2} <
        {{(m+3)(m+4)}\over 2}
$$
so $d<m+2$, and proposition \ref{horace} concludes.
\end{proof}

\begin{Lem}
\label{ultimpas}
If $i \leq m$, $d>m$, $4\,(i+j) \geq m^2-2\,m-4$ and $\m$ is
consistent in $U_s$ then $\ell_d(U_s,\m)=-1$.
\end{Lem}
\begin{proof}
By corollary \ref{ipetit} we may assume $2i >m$.
We will distinguish two cases according to the parity of $m$.
\begin{enumerate}
\item[$m$ even] Write $m=2\,t-2$. As $\m$ is consistent in $U_s$,
and $2i>m$, we must have $t \geq s$, so by semicontinuity
it is enough to see that there are no plane curves of
degree $d$ containing $Z_{K,\m}$ with $K$ general in
$U_{t}$.
It is clear that $(P_{t})_\m$ does not contain
$U_{t+1}$, and $4\,(i+j) \geq 4 \, (t^2 -3\, t+1)$,
therefore by proposition \ref{limit}
$$\V(t+1, m+1, i-t+1, j+t-2) \subset \V(t,m,i,j) \ .$$
By semicontinuity applied to the tautological flat
family on $\V(t,m,i,j)\subset \Hilb^{N}(S)$ then,
it is enough to see that there are no plane curves
of degree $d$ containing $Z_{K,\m'}$ with $K$ general
in $U_{t+1}$ and $\m'=(m+1,2^{i-t+1}, 1^{j+t-2})$.
But then $2(i-t+1)=2i-m \leq m$ and still
$2(i-t+1)+j+t-2 >m$ so the result follows by lemma \ref{foradosos}.

\item[$m$ odd] Write $m=2\,t-1$.
By semicontinuity it is enough to see that there are
no plane curves of degree $d$ containing $Z_{K,\m}$ with $K$ general
in $U_{t+1}$. But for $K \in U_{t+1}$, unloading
gives $Z_{K,\m}=Z_{K,(m+1,2^{i-t},1^{j+t})}$,
and again $2(i-t)=2i-m -1\leq m$ and 
$2(i-t)+j+t >m$ so the result follows by lemma \ref{foradosos}.
\end{enumerate}
\end{proof}

With this knowledge of level $d$ cluster schemes we are
in a position to attack the general case.
Let now $\rho$, $s$ be positive integers, and suppose we have
systems of multiplicities $\m_1, \m_2, \ldots, \m_\rho$ with
\begin{align*}
\m_1&=(m,2^{i_1},1^{j_1}) \\
\m_k&=(2^{i_k},1^{j_k}) \quad k=2, 3, \ldots, \rho .
\end{align*}
Suppose furthermore that $m \geq \min(s+\sum i_k,2s)$.
Consider a cluster scheme $Z_1=Z(K,\m_1)$ with $K \in U_s$ and
$\rho -1$ unibranched cluster schemes $Z_2, Z_3, \ldots,$ $Z_\rho$
supported at different points of $\P^2$, whose defining clusters
$K_2, K_3, \ldots, K_\rho$ have no satellite points.
The scheme $Z=Z_1 \cup Z_2 \cup \cdots \cup Z_\rho$ has length
$$ N={{m\,(m+1)} \over 2} + 3\, \sum_{k=0}^\rho i_k +
   \sum_{k=0}^\rho j_k \, .$$

\begin{Teo}
\label{rang}
If the position of the points of the $\rho$ clusters
$K_2, K_3, \ldots, K_\rho$ is general,
$3\, \sum i_k + \sum j_k \geq 2m+3$, and $4\, \sum j \geq m^2-4\,m-6$, 
then $Z$ has maximal rank, except in the cases
\begin{itemize}
 \item  $m=2$, $\sum i_k=4$, $j_k=0\ \forall k$, 
 \item  $m=4$, $\sum i_k=6$, $j_k=0\ \forall k$.
\end{itemize}
\end{Teo}
\begin {proof}
Consider
$\r=(r_1, r_2, \ldots, r_{\rho})$. Let $P_s(\r)$ be the closed
variety of $Y(\r)$ where $p_2, p_3, \ldots, p_s$ are proximate to
$p$. It is easy to see that $P_s(\r)$ is an irreducible variety
and $P_s=P_s(\r) \cap Y_{r_{\rho}-1}$. Then the claim may be
equivalently stated as
$$
\ell_d(P_s(\r),\m_0)= \max \left( -1, {{(d+1)(d+2)}\over 2}-N
   \right) \ \forall d
$$
where $\m_0=(m,2^{i_1},1^{j_1},2^{i_2},1^{j_2}, \ldots,
   2^{i_\rho},1^{j_\rho})$.
By semicontinuity it is enough to see that 
$\ell_d(P_s,\m_0)= \max ( -1, (d+1)(d+2)/2 -N ) \forall d$
and unloading gives that
$\m=(m,2^{\sum i_k},1^{\sum j_k})$ is equivalent to
$\m_0$ in $U_s$. Therefore
$$
\ell_d(P_s(\r),\m_0) \leq \ell_d(U_s,\m_0)=\ell_d(U_s,\m) \ \forall d \ .
$$
It only remains to be seen that
for a general cluster $K \in U_s$, $Z_{K,\m}$ has maximal
rank. Because of lemma \ref{nivelld} we can assume that $Z_{K,\m}$
has level $d>m$, and it is enough to see that 
$\ell_d(U_s,\m)= -1$.

Let $i=\sum i_k$, $j=\sum j_k$.
The proof runs by induction on $i$. The case $i \leq m$ has
already been settled in lemma \ref{ultimpas}, so suppose
$i \geq m+1$.
We will distinguish two cases according to the parity of $m$.
\begin{enumerate}
\item[$m$ even] Write $m=2\,t-2$. As $\m$ is consistent
in $U_s$ and $i \geq m+1$ we must have $t\geq s$, so
by semicontinuity it is enough to see that there are
no plane curves of degree $d$ containing $Z_{K,\m}$ with $K$ general
in $U_{t}$. It is clear that $(P_{t})_\m$ does not contain
$U_{t+1}$, and $4\,(i+j) \geq 4 \, (t^2 -3\, t+1)$,
therefore by proposition \ref{limit}
$$\V(t+1, m+2, i-2t+1, j+2t-2) \subset \V(t,m,i,j) \ .$$
By semicontinuity applied to the tautological flat
family on $\V(t,m,i,j)\subset \Hilb^{N}(S)$ then,
it is enough to see that there are no plane curves
of degree $d$ containing $Z_{K,\m'}$ with $K$ general
in $U_{t+1}$ and
$$\m'=(m', 2^{i'}, 1^{j'})=(m+2, 2^{i-2t+1}, 1^{j+2t-2}) \, .$$
An easy computation shows that we are still in the
numerical conditions of the claim, and $i'<i$, so we
can apply the induction hypothesis, except in the
case that $d=m+2$. In this case, either we are in
one of the exceptions above, or $i=m+1$, $j=3$, 
$\m'=(m+2,1^{m+3})$ and the claim follows from lemma
\ref{totuns}.

\item[$m$ odd] Write $m=2\,t-1$.
By semicontinuity it is enough to see that there are
no plane curves of degree $d$ containing $Z_{K,\m}$ with $K$ general
in $U_{t+1}$. But for $K \in U_{t+1}$, unloading
gives $Z_{K,\m}=Z_{K,(m+1, 2^{i-t}, 1^{j+t})}$,
so the result follows by the induction hypothesis again.
\end{enumerate}
\end{proof}

Remark that theorem \ref{rang} applies in particular
when we have only double points, proper or infinitely
near, and they are neither 5 nor less than 4. For the
remaining cases the behaviour is also known (cf. \cite{bh1}).
Namely, one double point has (obviously) maximal rank
and the scheme of three double points (proper or
infinitely near, in general position) has maximal rank;
systems which fail to have maximal rank appear in degree 2
when there are two points and in degree 4 when there are five.
In both cases there is an ``unexpected"
curve of the form $C=2D$, where $D$ is the straight line
or the conic through the points, respectively.
The cases of one point
of multiplicity $m=3, 4$ or 5 and $i$ double points
(again, proper or infinitely near) are also covered,
except for $i < 2m/3+1$ and for the case $m=4, i=6$.
All excepted cases involve less than eight points and were
therefore also solved by B. Harbourne. The systems which
do not have maximal rank are: $(3,2)$ in degree 3,
$(4,2)$, $(4,2^2)$ in degree 4, $(5,2)$,
$(5,2^2)$ in degree 5,
and $(4, 2^6)$  in degree 6.

We will now apply this result to find irreducible curves
of low degree with tacnodes and cusps. We will use a form of
Bertini's theorems slightly different from the usual ones.
Given a linear system $\L$ of curves on $S$ with no fixed
part, the base points of $\L$ form a (usually non unibranched)
cluster $\BP(\L)$.

\begin{Pro}
\label{bertini}
Let $\L$ be a linear system of curves on $S$ with no fixed
part. Then
\begin{enumerate}
  \item General curves in $\L$ go sharply through $\BP(\L)$.
  \item If $\L$ is reducible then it is composed of the curves
      of a pencil.
\end{enumerate}
\end{Pro}
\begin{proof}
2 and the fact that all singularities of general curves sit
at the proper base points of $\L$ are standard. For a proof
of 1, cf. \cite[7.2]{cll}. Cf. also
Zariski's remark on the theorems of Bertini in \cite{zas}.
\end{proof}

\begin{Lem}
\label{nobase}
Let $Z \subset Z'$ be two zero-dimensional subschemes of $\P^2$
such that $\len Z' = \len Z +1$, and let $d$ be a positive integer
such that going through $Z$ imposes independent conditions to
curves of degree $d$. Then $Z'$ imposes independent conditions
to curves of degree $d+1$.
\end{Lem}
\begin{proof}
This is an easy application of the residual exact
sequence of the Horace method.
There is a unique point $p \in \P^2$ where the length of the component
of $Z'$ supported at $p$ is bigger than that of $Z$. Let $L$ be
a general straight line through $p$, and consider the residual
exact sequence
$$
0 \longrightarrow \I_{Z''}(d) \longrightarrow \I_{Z'}(d+1)
  \longrightarrow \I_{Z' \cap L / L}(d+1) \longrightarrow 0
$$
We have to prove that
$$H^1(\I_{Z' \cap L / L}(d+1)) = H^1(\I_{Z''} (d)) = 0 \ .$$
As going through $Z$ imposes independent conditions to curves of
degree $d$ and $L$ is general, $\len Z \cap L \leq d+1$, therefore
$\len Z' \cap L \leq d+2$ and $H^1(\I_{Z' \cap L / L}(d+1))=0$.
As by hypothesis $H^1(\I_Z (d))=0$, it will be enough to see
that $Z$ contains the residual scheme $Z''$. Let
$I, I', I'' \subset \O_p$ be the ideals locally defining
$Z$, $Z'$ and $Z''$, and let $f \in \O_p$ be a local equation
of $L$. $\len Z' = \len Z +1$ implies that $I=I'+(g)$ for
some $g \in \O_p$ with $g \M_p \subset I'$, therefore
$ f I = fI' + (fg) \subset I'$ and $I \subset I'' = (I':f)$,
as wanted.
\end{proof}

\begin{Cor}
\label{bertiniclusters}
If $(K,\m)$ is a weighted (not necessarily unibranched) cluster 
such that $Z_{K,\m}$ imposes independent
conditions to curves of degree $d$, then general curves in
$\L=\L_{d+1}(K,\m)$ go sharply through $(K,\m)$. Furthermore,
if $(K,\m)$ is not a single point of multiplicity $d+1$
then a general curve in $\L$ is irreducible.
\end{Cor}
\begin{proof}
Because of lemma \ref{nobase}, $\L$ has no fixed
part and $\BP(\L)=(K,\m)$. By proposition \ref{bertini}
then, we only have to see that if $\L$ is composed
of the curves of a pencil then $\BP(\L)$ is a point with multiplicity
$d+1$. Let $\L$ be composed of $r$ curves of degree
$k$ in a pencil. Then $rk=d+1$ and
$$
r=\ell_{d+1}(K,\m)=d+2 + \ell_{d}(K,\m)\geq d+1
$$
because $Z_{K,\m}$ imposes independent conditions to
curves of degree $d$. Therefore $r=d+1$ and $k=1$.
A pencil of lines has a unique base point $p$
so curves in $\L$ are composed of $d+1$ lines through
$p$ and $\BP(\L)$ is as claimed.
\end{proof}

\begin{figure}
  \begin{center}
    \mbox{\includegraphics[2122,2444][2430,2674]{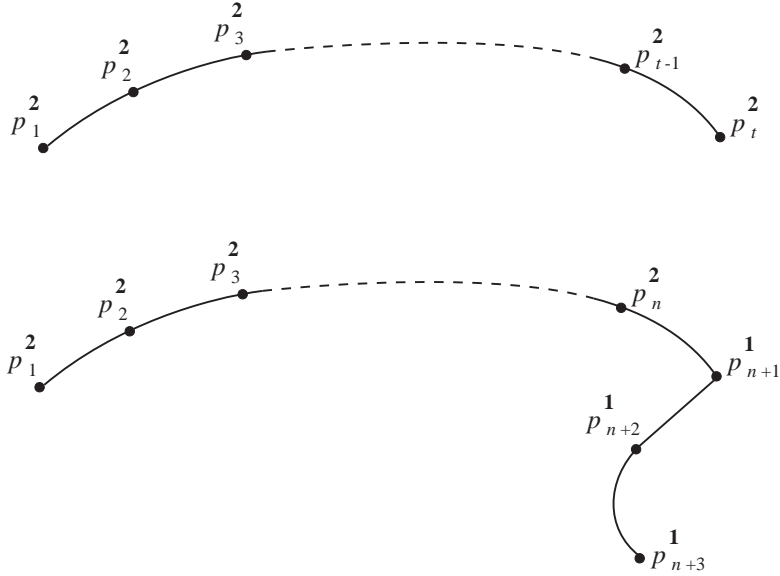}}
    \caption{Enriques diagrams of a tacnode cluster
         and an extended cusp cluster.}
  \end{center}
\end{figure}

The cluster of infinitely near singular points of a tacnode
of order $t$ is a weighted unibranched cluster $K$ of $t$ free
points with multiplicities $\m=(2^t)$. The scheme $Z_{K,\m}$
is called a tacnode scheme.
The cluster of infinitely near singular points of a cusp
of order $n$ is a weighted unibranched cluster
$K$ of $n+2$ points, the last of which being satellite
and the others free, with multiplicities
$\m=(2^n,1^2)$. We define an ``extended" unibranched
cluster $(K',\m')$ which has an additional free point taken
with multiplicity one and call $Z_{K',\m'}$ a cusp scheme of
order $n$. When this last point $p_{n+3}\in K'$ varies in
$E_{n+2} \cong \P^1$ we obtain a flat family of cluster
schemes, because of lemma \ref{unpuntmes}.
In the special position of $p_{n+3}$ which makes
it proximate to $p_{n+1}$, unloading gives
$\delta(\m')=(2^{n+1},0,0)$ so any cusp scheme
of order $n$ can be specialized to 
the tacnode scheme of order $n+1$ given by the
free points $p_1, \ldots, p_{n+1}$.

\begin{Teo}
\label{existencia}
If
$$ {{d(d+1)} \over 6} \geq  \sum_{i=1}^\tau t_i +
    \sum_{i=1}^\nu (n_i +1) \geq 6$$
then there exists a reduced irreducible curve of
degree $d$ with $\tau$ tacnodes of orders $t_1, \ldots, t_\tau$
and $\nu$ cusps of orders $n_1, \ldots, n_\nu$ as
its only singularities.
\end{Teo}

\begin{proof}
Consider a scheme
$$ Z=T_1 \cup T_2 \cup \cdots \cup T_\tau \cup
     N_1 \cup N_2 \cup \cdots \cup N_\nu $$
where $T_i$ is a tacnode scheme of order $t_i$ supported
at $p_i \in \P^2$ and $N_i$ is a cusp scheme of order $n_i$
supported at $q_i$, all of them having their points in
general position.
We claim that the linear system $\L_d(Z)$
of curves of degree $d$ containing $Z$ is nonempty and that
a generic curve in it has a tacnode of order $t_i$ at $p_i$,
a cusp of order $n_i$ at $q_i$ and no other singularities.

Specializing the cusp schemes $N_i$ to tacnode schemes
$\bar T_i$ of order $n_i+1$ we obtain a scheme $\bar Z$ in the
conditions of theorem \ref{rang}, therefore of maximal
rank and length $3 \left( \sum t_i + \sum (n_i +1) \right)$.
So by semicontinuity $Z$ is also of maximal rank. The
bound on $d$ assures that $Z$ imposes independent conditions
to curves of degree $d-1$. So
by \ref{bertiniclusters} a general curve in $\L_d(Z)$ is irreducible
and has no other singularities but the ones in $p_i$, $q_i$, which
are tacnodes and cusps of the desired orders.
\end{proof}

Remark that in the hypotheses of theorem \ref{existencia}
we ask that $\sum t_i +\sum (n_i +1) \geq 6$
in order to apply \ref{rang}. In fact, the
remarks we made after the proof of theorem \ref{rang}
prove that the result holds also for
$\sum t_i +\sum (n_i +1) =3, 4$. In case
$\sum t_i +\sum (n_i +1) =5$, an ad--hoc reasoning
can be used to prove the existence; however, for
these small numbers of singularities the result is
neither new nor significant, so we omit this.

We would like to point out that there are examples in
\cite{los} of curves $F_k$ with one tacnode or cusp, of degree
lower than the one given by \ref{existencia}. Namely,
they have degree $d=2k +1$ and
$$\begin{cases}
\text{a tacnode of order } \frac{2k^2+3k-1}{2} \text{ if }
   2k^2+3k-1 \text{ is even} \\
\text{a cusp of order } \frac{2k^2+3k-2}{2} \text{ if }
   2k^2+3k-1 \text{ is odd} \\
\end{cases}$$
The author was informed by C. Lossen that he conjectures,
after testing many particular cases, that these curves are
irreducible and have no other singularity. He conjectures
also that they do \emph{not} satisfy the \emph{T--smoothness}
property (cf. \cite{shu}) which states that
the variety of curves of degree $d$ with that singularity
is smooth of the expected dimension at $F$.
Curves with tacnodes and cusps whose points are in general
position (as those given by \ref{existencia})
do satisfy the T--smoothness property.

The reader may notice that after \ref{bertiniclusters},
any $h^1$--vanishing result for a class of cluster schemes
can be exploited to obtain curves of low degree with the
equisingularity type fixed by the clusters. In fact, the
bounds obtained using \ref{bertiniclusters} are (slightly)
sharper than those obtained using the somewhat more
complicated reasoning of \cite{gls}.
In particular, it is not difficult to extend theorem
\ref{existencia} to curves with tacnodes, cusps, and
one different singularity.
\begin{Pro}
Let $S$ be an equisingularity class whose cluster of
infinitely near singular points is unibranched and consists
of one point of multiplicity $m$ followed by $k$ free
double points. Define
$M= k + \sum_{i=1}^\tau t_i + \sum_{i=1}^\nu (n_i +1) 
  + \max \left( 0, {{m^2 - 4m - 6} \over 3}\right)$.
If
$$
 {{d(d+1)} \over 6} - {{m(m+1)} \over 6} \geq M \ge {2 \over 3} m+1
$$
then there exists a reduced irreducible curve of
degree $d$ with one singularity of type $S$,
$\tau$ tacnodes of orders $t_1, \ldots, t_\tau$
and $\nu$ cusps of orders $n_1, \ldots, n_\nu$ as
its only singularities.
\end{Pro}
\begin{proof}
Analogous to the proof of \ref{existencia}.
\end{proof}

Remark also that in the proof of theorem \ref{existencia} we proved
that the union of general tacnode and cusp schemes has
maximal rank, which is a case not included in theorem
\ref{rang}. It is easy to see that many other cases may
be equally treated using the same techniques. For example,

\begin{Pro}
The scheme of the cluster of infinitely
near singular points of a $D_k$--singularity
whose points are in general position has maximal rank,
except in the two cases $k=6,7$.
\end{Pro}
\begin{proof}
The cluster of infinitely near singular points of a $D_k$--singularity,
with $k$ even, is a weighted unibranched cluster with
$k/2-1$ free points with multiplicities $\m=(3,2^{k/2-2})$,
so we are in a particular case of theorem \ref{rang}.

If $k$ is odd, then the cluster of infinitely near singular points
of $D_k$ is a weighted unibranched cluster with $r=(k+1)/2$ points,
the last of which being satellite and the others free (that is,
$K \in U_{r,r-1}$), with multiplicities $\m=(3,2^{r-3},1^2)$.
For simplicity, we assume $k\ge 13$ (or $r \ge 7$), as the cases
with $k$ small require special care.
Let $U \subset P_{r,r-1}\cap P_3$ be the open subset where
only $p_3$ and $p_r$ are satellites (we specialize the
third point to be proximate to the first). For $K\in U$,
unloading gives $Z_{K,\m}=Z_{K,\m'}$ with $\m'=(4,2^{r-4},1,0,0)$,
so the length of $Z_{K,\m}$ is the same for $K$ in $U$
or in $U_{r,r-1}$, therefore $U \subset (P_{r,r+1})_\m$.
As $U_{r,r-1}$ is dense in $U \subset (P_{r,r+1})_\m$,
by semicontinuity it is enough to see that
$Z_{K,\m}=Z_{K,\m'}$ has maximal rank for $K$
general in $U$
and the result follows from theorem \ref{rang},
because the multiplicity of the last point
is now 0 so we can assume $K$ general in $U_3$.
\end{proof}


\begin{thebibliography}{}
   \bibitem{ah}    Alexander, J. Hirschowitz, A.,
                        Interpolation on Jets
                        \emph{J. Alg.},
                        \textbf{192} (1997), 412-417.
   \bibitem{bar}    Barkats, D.,
                        Vari\'et\'es des courbes planes \`a noeuds
                        et \`a cusps, in: Peter E. Newstead, ed.,
                        \emph{Algebraic Geometry,} 
                        Lect. Notes Pure App. Math. \textbf{200},
                        M. Dekker 1988, pp. 25-36.
   \bibitem{bri}    Brian‡on, J.,
                        Description de $\Hilb^n \C\{x,y\}$,  
                        \emph{Inventiones Math.}
                        \textbf{41} (1977), 45-89.
   \bibitem{cma}    Casas-Alvero, E.,           
                        Infinitely near imposed singularities and 
                        singularities of polar curves,
                        \emph{Math. Ann.}
                        \textbf{287} (1990), 429-454.
   \bibitem{cll}     Casas-Alvero, E.,
                        \emph{Singularities of plane curves}
                        preprint Universitat de Barcelona (1997).
   \bibitem{cat}    Catalisano, M.V., Gimigliano, A.
                        On curvilinear subschemes of $\P^2$ 
                        \emph{J. Pure App. Algebra}
                        \textbf{93} (1994), 1-14.
   \bibitem{cm1}     Ciliberto, C. Miranda, R.,
                        Interpolation on Curvilinear Schemes,
                        \emph{J. Alg.}
                        \textbf{203} (1998) 677-678.
   \bibitem{cm2}     Ciliberto, C. Miranda, R.,
                        Degenerations of planar linear systems,
                        \emph{Journal Reine Ang. Math.}
                        \textbf{501} (1998) 191-220.
   \bibitem{cm3}     Ciliberto, C. Miranda, R.,
                        Linear systems of plane curves with base
                        points of equal multiplicity,
                        (1998) to appear in \emph{Trans. A. M. S.}
   \bibitem{enr}    Enriques, F. - Chisini, O.,
                        \emph{Lezioni sulla teoria geometrica
                        delle equazioni e delle funzioni algebriche,}
                        N. Zanichelli, Bologna 1915.
   \bibitem{ev1}    Evain, L.
                        \emph{Collisions de trois gros points sur
                        une surface alg\'ebrique}, thesis, Nice 1997. 
   \bibitem{ev3}      Evain, L.,
                        Une minoration du degr' des courbes planes
                        … singularit's impos'es,
                        \emph{Preprint ENS Lyon,}
                        \textbf{212} (1997), 1-17.
   \bibitem{fog}     Fogarty, J.,
                        Algebraic families on an     
                        algebraic surface,
                        \emph{Amer. J. Math}
                        \textbf{10} (1968) 511-521.
   \bibitem{gls}   Greuel, G.M. - Lossen, C. - Shustin, E.,
                        Plane curves of minimal degree 
                        with prescribed singularities,
                        \emph{Inventiones Math.}                        
                        \textbf{133} (1998), 539-580.
   \bibitem{gls2}   Greuel, G.M. - Lossen, C. - Shustin, E.,
                        Castelnuovo function, zero-dimensional
                        schemes and singular plane curves,
                        \emph{duke e--print}
                        \textbf{9903179} (1999).
   \bibitem{ega}        Grothendieck, A. - Dieudonn', J.,
                        EGA IV, 4, 
                        \emph{El'ments de g'om'trie alg'brique,}
                        Inst. Hautes Etudes Sci. Publ. Math.
                        \textbf{32} (1967).
   \bibitem{bh1}         Harbourne, B.,
                        Complete linear systems on rational
                        surfaces,
                        \emph{Trans. A.M.S.}
                        \textbf{289} (1985) 213-226.
   \bibitem{bh2}         Harbourne, B.,
                        The geometry of rational surfaces
                        and Hilbert functions of points in
                        the plane,
                        \emph{Can. Math. Soc. Conf. Proc.} 
                        vol. \textbf{6} (1986) 95-111.
   \bibitem{hasun}    Harbourne, B.,
                        Iterated blow-ups and moduli    
                        for rational surfaces, in:
                        A. Holme and R. Speiser, eds.,
                        \emph{Algebraic Geometry Sundance 1986},
                        LNM \textbf{1311}, Springer 1988, pp. 101-117.
   \bibitem{bh3}     Harbourne, B.,
                        Rational surfaces with $K^2>0$,
                        \emph{Proc. A. M. S.}
                        \textbf{124} (1996) 727-733.
   \bibitem{bh4}         Harbourne, B.,
                        Anticanonical rational surfaces,
                        \emph{Trans. A.M.S.}
                        \textbf{349} (1997) 1191-1208.
   \bibitem{hag}   Hartshorne, B.,
                        \emph{Algebraic Geometry},
                        GTM 52, Springer 1977.
   \bibitem{hir}    Hirschowitz, A.,
                        La M'thode d'Horace pour l'interpolation
                        a plusieurs variables,
                        \emph{Manuscripta Math.}
                        \textbf{50} (1985), 337-388.
   \bibitem{hir2}    Hirschowitz, A.,
                        Une conjecture pour la cohomologie
                        des diviseurs sur les surfaces
                        rationnelles generiques,
                        \emph{J. Reine Angew. Math}
                        \textbf{397} (1989), 208-213.
   \bibitem{iar}     Iarrobino, A.
                        Punctual Hilbert Schemes,
                        \emph{Mem. Amer. Math. Soc.}
                        \textbf{188} (1977).
   \bibitem{klit}    Kleiman, S.L.,
                        Multiple-point formulas I: Iteration,
                        \emph{Acta Math.}
                        \textbf{147} (1981), 13-49.
   \bibitem{los}    Lossen, C.,
                        New asymptotics for the existence of plane
                        curves with prescribed singularities,
                        to appear in \emph{Comm. in Algebra} (1999).
   \bibitem{mir}  Miranda, R.,
                        Linear systems of plane curves,
                        \emph{Notices of the A. M. S.}
                        \textbf{46}, 2, (1999), 192-202.
   \bibitem{nv} Nobile, A., Villamayor, O.,
                        Equisingular stratifications associated
                        to families of planar ideals,
                        \emph{J. Alg.}
                        \textbf{193} (1997), 239-259.
   \bibitem{pax} Paxia, G.,
                        On flat families of fat points,
                        \emph{Proc. A. M. S.}
                        \textbf{112} (1991), 19-23.
   \bibitem{ran}   Ran, Z.,
                        Curvilinear enumerative geometry,
                        \emph{Acta Math.}
                        \textbf{155} (1985), 81-101.
   \bibitem{ran1}   Ran, Z.,
                        Enumerative geometry of singular
                        plane curves,
                        \emph{Inventiones Math.}
                        \textbf{97} (1989), 447-469.
   \bibitem{ran2}   Ran, Z.,
                        On the Nagata problem,
                        \emph{duke e--print}
                        \textbf{9809101} (1998).
   \bibitem{meu}  Ro\'e, J.,
                        On the existence of plane curves with
                        imposed multiple points,
                        \emph{duke e--print}
                        \textbf{9807066} (1998).
   \bibitem{shu}  Shustin, E.,
                        Smoothness of equisingular families
                        of plane algebraic curves,
                        \emph{Int. Math. Res. Notices}
                        (1997), 67-82.
   \bibitem{zas}   Zariski, O.,
                        \emph{Algebraic Surfaces},
                        2nd. suppl. ed., Ergebnisse 61,
                        Springer 1971.
   \bibitem{zpm}   Zariski, O.,
                        \emph{Le problŠme des modules
                        pour les branches planes},
                        Hermann 1986.
\end{thebibliography}
\end{document}